\documentclass{article}

\usepackage{amsthm}
\usepackage{amsfonts}
\usepackage{amssymb}
\usepackage{amsgen}
\usepackage{amsmath}
\usepackage{amsopn}
\usepackage{verbatim}
\usepackage{xypic}
\usepackage{xspace}
\usepackage{multicol}
\usepackage{makeidx}
\usepackage{eepic}
\usepackage{upref}

\theoremstyle{plain}
\newtheorem{thm}{Theorem}[section]
\newtheorem{lem}[thm]{Lemma}
\newtheorem{prop}[thm]{Proposition}
\newtheorem{cor}[thm]{Corollary}
\newtheorem{conj}[thm]{Conjecture}
\theoremstyle{definition}

\newtheorem{defn}[thm]{Definition}
\newtheorem{rem}[thm]{Remark}

\newtheorem{mainconj}{Main Conjecture}

 \font\cyr=wncyr10
 \newcommand{\nc}{\newcommand}
 \nc{\reg}{{\rm reg}}
 \nc{\mal}{{{\scriptstyle \maltese}}}
 \nc{\re}{{\Re}}
 \nc{\fA}{{\mathfrak A}}
 \nc{\ra}{\rightarrow}
 \nc{\ors}{{\vec s\,}}
 \nc{\os}{{\overset}}
 \nc{\Z}{{\mathbb Z}}
 \nc{\R}{{\mathbb R}}
 \nc{\N}{{\mathbb N}}
 \nc{\ZN}{{\mathbb Z_{\ge 0}}}
 \nc{\Q}{{\mathbb Q}}
 \nc{\C}{{\mathbb C}}
 \nc{\D}{{\mathcal D}}
 \nc{\caT}{{\mathcal T}}
 \nc{\tB}{{\tilde B}}
 \nc{\Li}{{\rm Li}}
 \nc{\suf}{{\ast\,}}
 \nc{\sufq}{{\ast_q\,}}
 \nc{\gam}{{\gamma}}
 \nc{\ga}{{\alpha}}
 \nc{\gl}{{\lambda}}
 \nc{\gb}{{\beta}}
 \nc{\gd}{{\delta}}
 \nc{\gs}{{\sigma}}
 \nc{\gS}{{\Sigma}}
 \nc{\sif}{{\mathcal S}}
 \nc{\gt}{{\tau}}
 \nc{\Lra}{\Longrightarrow}
 \nc{\lra}{\longrightarrow}
 \nc{\fS}{{\mathfrak S}}
 \nc{\DD}{{\mathfrak D}}
 \nc{\Llra}{\Longleftrightarrow}
 \nc{\ol}{\overline}
 \nc{\lms}{\longmapsto}
 \nc{\zq}{{\zeta_q}}
 \nc\qup{{q\uparrow 1}}
 \nc{\us}{\underset}
 \nc{\tn}{{\tilde{n}}}
 \nc{\gD}{{\Delta}}
 \nc{\bk}{{\bf k}}

\nc{\sha}{{\mbox{\cyr x}}}

\begin{document}

\title{Double Shuffle Relations of Euler Sums}
\author{Jianqiang Zhao}
\date{}
\maketitle
\begin{center}
{\large Department of Mathematics, Eckerd College, St. Petersburg,
FL 33711}
\end{center}
\vskip0.6cm

\noindent{\bf Abstract.} In this paper we shall develop a theory
of (extended) double shuffle relations of Euler sums which
generalizes that of multiple zeta values (see Ihara, Kaneko and
Zagier, \emph{Derivation and double shuffle relations for multiple
zeta values}.  Compos.~Math.~\textbf{142} (2)(2006), 307--338).
After setting up the general framework we provide some numerical
evidence for our two main conjectures. At the end we shall prove
the following long standing conjecture: for every positive integer
$n$
$$\zeta(\{3\}^n)=8^n\zeta(\{\ol2,1\}^n).$$
The main idea is to use the double shuffle relations and the
distribution relation. This particular distribution relation
doesn't follow from the double shuffle relation in general. But we
believe it does follow from the extended double shuffle relations.

\section{Introduction}
There are many different generalizations of Riemann zeta
functions. One may introduce more variables to define the multiple
zeta function as
\begin{equation}\label{zeta}
\zeta(s_1,\dots, s_l)=\sum_{k_1>\dots>k_l>0}
 \frac{1}{k_1^{s_1}\cdots k_l^{s_l}}
\end{equation}
for complex variables $s_1,\dots, s_l$ satisfying
$\re(s_1)+\dots+\re(s_j)>j$ for all $j=1,\dots,l$.  It was Euler
who first systematically studied the special values of these
functions at positive integers when $d=2$, after corresponding
with Goldbach. Among many results he showed (see \cite{LE1} and
\cite[p.\ 266]{LE2}),
 $$2\zeta(m,1)=m\zeta(m+1)-\sum_{j=1}^{m-2}\zeta(j+1)\zeta(m-j),
    \qquad 2\le m\in\Z.$$
However, only in the past fifteen years or so have these values
been found to have significant arithmetic, algebraic and geometric
meanings and have since been under intensive investigation (see
\cite{MG,H1,LM,Zag}). Consequently many other multiple zeta value
(MZV) identity families have been discovered and it is conjectured
\cite{IKZ} that all of them are consequences of the finite and
extended double shuffle relations (see section \ref{sec:EDS} for
details).

In another direction, MZVs can also be thought of as special
values of the multiple polylogarithms (note that $s_i$ are all
positive integers and $s_1>1$)
\begin{equation}\label{polylog}
Li_{s_1,\dots, s_l}(x_1,x_2,\dots,x_l)=\sum_{k_1>\dots>k_l>0}
\frac{x_1^{k_1}\cdots x_l^{k_l}}{ k_1^{s_1}\cdots k_l^{s_l}}.
\end{equation}
Goncharov \cite{Gcyc} proposes to study the special values of
these functions at roots of unity and believes this will provide
the high cyclotomic theory. Moreover, theoretical physicists have
already found out that such values appear naturally in the study
of Feynmen diagrams (\cite{Br1,Br2}). We will study these special
values in another paper \cite{Zh}.

Starting from early 1990's Hoffman \cite{H1,H2} has constructed
some quasi-shuffle\footnote{We will call ``stuffle'' in this
paper.} algebras in order to catch the essence of MZVs. Recently
he \cite{H3} extends this to incorporate the special values of
polylogarithms at roots of unity, although his definition of
$*$-product is different from ours. If we only take $x_i=\pm 1$ in
the multiple polylogarithms then the special values
$Li_{s_1,\dots, s_l}(x_1,x_2,\dots,x_l)$ are called {\em
(alternating) Euler sums} (see \cite{BBB1}):
\begin{equation}\label{equ:z}
\zeta(s_1,\dots,s_l;x_1,\dots,x_l)
   := \sum_{k_1>\cdots>k_l>0}\;\prod_{j=1}^l
   \frac{x_j^{k_j}}{k_j^{s_j}}.
\end{equation}
We will only consider such sums in this paper. Observe that we may
even allow $s_1=1$ if $x_1=-1$. To save space, if $x_j=-1$ then
$\overline s_j$ will be used and if a substring $S$ repeats $n$
times in the list then $\{S\}^n$ will be used. For example,
$\zeta(\overline1)=\zeta(1;-1)=-\ln 2$ and $\zeta(2)=\pi^2/6$. We
will call indices like $(\bar 1,2,\bar3)$ \emph{signed indices.}

It is well known that there are two types of relations among MVZs,
one from multiplying the series \eqref{equ:z} and the other from
multiplying their iterated integral representations. Both of these
can be generalized to Euler sums fairly easily. After briefly
sketching this theory in section \ref{sec:EDS} and posing two
conjectures we shall provide some numerical computation to support
them in section \ref{sec:comp}.

The rest of the paper is devoted to the proof of
\begin{thm} {\em For every positive integer} $n$
$$\zeta(\{3\}^n)=8^n\zeta(\{\ol2,1\}^n).$$
\end{thm}
Around 1996 Borwein, Bradley and Broadhurst~\cite{BBB2} first
noticed that the above result must be true after some intensive
computation. It is remarkable that this was the only conjectured
family of identities relating alternating Euler sums to MZVs.
Several proofs of the case $n=1$ can be found in \cite{BdB}. The
case $n=2$ is much more difficult and the only known proof before
this work was by computer computation \cite{B}. In this paper, we
will prove this result in general by using double shuffle
relations and the distribution relation. However, in general it is
impossible to prove the identities by just the finite double
shuffle relations.

I would like to thank David Bradley for his encouragement and many
email discussions. In particular, he pointed out the equivalent
form of Theorem \ref{thm:iden} in Theorem \ref{thm:key}. This
simplifies my original computation greatly.

\section{The double relations and the algebra $\fA$}\label{sec:EDS}
Kontsevich first noticed that MZVs can be represented by iterated
integrals. It is quite natural and easy to extend this to Euler
sums (see \cite{BdB}). Set
 $$a=\frac{dt}{t},\qquad b=\frac{dt}{1-t},\qquad c=\frac{-dt}{1+t}.$$
For every positive integra $n$ define
$$\gb_n=a^{n-1}b \quad\text{and}\quad \gam_n=a^{n-1}c.$$
Then it is straight-forward to verify that for
$s_1>1$\begin{equation}\label{equ:mzv}
\zeta(s_1,\dots,s_l)=\int_0^1 \gb_{s_1}\cdots \gb_{s_l}:=\int_0^1
\gb_{s_1}(t_1) \left(\int_0^{t_1} \gb_{s_2}(t_2)\cdots
\int_0^{t_{l-1}} \gb_{s_l}(t_l) d\, t_l \cdots d\, t_2 \right)
d\,t_1
\end{equation}
To study this for general Euler sums we can follow Hoffman
\cite{H2} by defining an algebra of words as follows:
\begin{defn} Set $A_0=\{\bf 1\}$ to be the set of empty
word. Define $\fA=\Q\langle A\rangle$ to be the graded
noncommutative polynomial $\Q$-algebra generated by letters $a,$
$b$ and $c$, where $A$ is a locally finite set of generators whose
degree $n$ part $A_n$ consists of words (i.e., a monomial in the
letters) of length $n$. Let $\fA^0$ be the subalgebra of $\fA$
generated by words not beginning with $b$ and not ending with $a$.
The words in $\fA^0$ are called \emph{admissible words.}
\end{defn}

Observe that every Euler sum can be expressed as an iterated
integral over $[0,1]$ of a unique admissible word $w$ in $\fA^0$.
Then we denote this Euler sum by $Z(w)$. It is quite easy to see
that $\fA^0$ is generated by words $\gb_n$ ($n\ge 2$) and $\gam_m$
($m\ge 1$). For example from \eqref{equ:mzv}
\begin{align*}
\zeta(s_1,\dots,s_l)=Z(\gb_{s_1}\cdots \gb_{s_l})
\end{align*}
If some $s_i$'s are replaced by $\bar s_i$'s then we need to
change some $\gb$'s to $\gam$'s according to the following:

\begin{quote}
{\bf Converting rule between signed indices and admissible words
in $\fA^0$}. Going down from $s_1$ to $s_l$, as soon as we see the
first signed letter $\bar s_i$ we change every $\gb$ after
$\gb_{s_i}$ (inclusive) to $\gam$ until the next signed letter
$\bar s_j$ is encountered. We then leave alone and all the $\gb$'s
after $\gb_{s_j}$ (again inclusive) until we see the next signed
letter when we start to toggle again. Carry on this toggling till
the end.
\end{quote}
Imaginatively we can think the bars as switches between $\gam$'s
and $\gb$'s. It is not hard to see that this establishes a
one-to-one correspondence between Euler sums and the words in
$\fA^0$. For example:
\begin{align*}
\zeta(\bar1,2,2,\bar4,3,\bar5,\bar6)
=Z(\gam_1\gam_2\gam_2\gb_4\gb_3\gam_5\gb_6)=Z(cacaca^3ba^2ba^4ca^5b).
\end{align*}
We would like to find many relations between different special
values. Remarkably, Chen \cite{Chen} developed a theory of
iterated integral which can be applied in our situation.
\begin{lem}\label{chen's}
Let $w_i$ $(i\ge 1)$ be $\C$-valued 1-forms on a manifold $M$. For
every path $p$,
 $$ \int_p w_1\cdots w_r\int_p w_{r+1}\cdots w_{r+s}=
 \int_p (w_1\cdots w_r) \sha (w_{r+1}\cdots w_{r+s})$$
where $\sha$ is the shuffle product defined by
 $$(w_1\cdots w_r) \sha (w_{r+1}\cdots
 w_{r+s})=\sum_{\substack{\gs\in
 S_{r+s},\gs^{-1}(1)<\cdots<\gs^{-1}(r)\\
 \gs^{-1}(r+1)<\cdots< \gs^{-1}(r+s)}}
   w_{\gs(1)}\cdots w_{\gs(r+s)}.$$
\end{lem}
For example, we have
 $$\zeta(\bar1)\zeta(2)=Z(c)Z(ab)=Z(c\sha (ab))=Z(cab+acb+abc)
 =\zeta(\bar1,\bar2)+\zeta(\bar2,\bar1)+\zeta(2,\bar1).$$
Let $\fA_\sha$ be the algebra of $\fA$ together with the
multiplication defined by shuffle product $\sha$. Denote the
subalgebra $\fA^0$ by $\fA_\sha^0$ when we consider the shuffle
product. Then we can easily prove
\begin{prop} \label{shahomo} The map $Z:  \fA_\sha^0\lra \R$,
is an algebra homomorphism.
\end{prop}

On the other hand, it is well known that Euler sums also satisfy
the series stuffle relations. For example
\begin{equation*}
   \zeta(\bar1)\zeta(2)=\zeta(\bar1,2)+\zeta(2,\bar1)+\zeta(\bar3).
\end{equation*}
because
$$\sum_{j>0}\sum_{k>0}=\sum_{j>k>0}+\sum_{k>j>0}+\sum_{j=k>0}.$$
To study such relations in general we need the following
definition.
\begin{defn}  Denote by $\fA^1$ the subalgebra of
$\fA$ which is generated by words $\gb_k$ and $\gam_k$ with $k\ge
1$. In other words, $\fA^1$ is the subalgebra of $\fA$ generated
by words not ending with $a$. For any word $w\in\fA^1$ and
positive integer $n$ define the maltese operator
$\mal_{\gb_n}(w)=w$, and $\mal_{\gam_n}(w)$ to be the word with
$\gb$ and $\gam$ toggled. For example
$\mal_{\gam_1}(\gam_2\gb_4)=\gb_2\gam_4.$ We then define a new
multiplication $*$ on $\fA^1$ by requiring that $*$ distribute
over addition, that $1*w=w*1=w$ for any word $w$, and that, for
any words $w_1,w_2$ and letters $x$ and $y$,
 \begin{equation}\label{equ:defnstuffle}
xw_1*yw_2 = x\Big(\mal_x\big(\mal_x(w_1)*yw_2\big)\Big)
 + y\Big(\mal_y\big(xw_1*\mal_y(w_2)\big)\Big) +
 [x,y]\Big(\mal_{[x,y]}\big(\mal_x(w_1)*\mal_y(w_2)\big)\Big)
\end{equation}
where
 $$[\gb_m,\gb_n]=[\gam_m,\gam_n]=\gb_{m+n},\quad
  [\gam_m,\gb_n]=[\gb_m,\gam_n]=\gam_{m+n}.$$
We call this multiplication the \emph{stuffle product}.
\end{defn}
If we denote $\fA^1$ together with this product $*$ by $\fA_*^1$
then it is not hard to show that
\begin{thm} \em{(Compare \cite[Theorem 2.1]{H2})}
The polynomial algebra $\fA_*^1$ is a commutative graded
$\Q$-algebra.
\end{thm}

Now we can define the subalgebra $\fA_*^0$ similarly to
$\fA_\sha^0$ by replacing the shuffle product by stuffle product.
Then by induction on the lengths and using the series definition
we can quickly check that for any $w_1,w_2\in \fA_*^0$
$$Z(w_1)Z(w_2)=Z(w_1\ast w_2).$$
This implies that
\begin{prop} \label{*homo}
The map $Z:  \fA_*^0  \lra \R$, is an algebra homomorphism.
\end{prop}

For $w_1,w_2\in \fA^0$ we will say that
 $$Z(w_1\sha w_2-w_1*w_2)=0$$
is a finite double shuffle (FDS) relation. It is known that even
for MZVs these relations are not enough to recover all the
relations among MZVs. However, we believe one can remedy this by
considering extended double shuffle relations produced by the
following mechanism. This wss explained very well in \cite{IKZ}
when Ihara, Kaneko and Zagier considered MZVs. So we will follow
them closely in the rest of the section.

Combining Propositions \ref{*homo} and \ref{shahomo} we can prove
easily (see \cite[\S2 Prop. 1]{IKZ}):
\begin{prop} \label{prop:eDS}
We have two algebra homomorphisms:
$$Z^*: (\fA_*^1,*)\lra \R[T],\quad \text{and}\quad Z^\sha: (\fA_\sha^1,\sha)\lra \R[T]$$
which are uniquely determined by the properties that they both
extend the evaluation map $Z:\fA^0\lra \R$ and send $b$ to $T$.
\end{prop}

For any signed index $\bk=(k_1,\dots,k_n)$ where $k_i$ are
positive integers (it may have a bar on top), let the image of the
corresponding words in $\fA^1$ under $Z^*$ and $Z^\sha$ be denoted
by $Z_\bk^*(T)$ and  $Z_\bk^\sha(T)$ respectively. For example,
$$\zeta(\bar1)T=Z_{\bar1}^*(T)Z_{1}^*(T)=Z^*(c*b)
 =Z_{(1,\bar1)}^*(T)+\zeta(\bar1,1)+\zeta(\bar2)$$
while
$$\zeta(\bar1)T=Z_{\bar1}^\sha(T)Z_{1}^\sha(T)
 =Z^\sha(c\sha b)=Z_{(1,\bar1)}^\sha(T)+\zeta(\bar1,\bar1).$$
From this and more computations we believe that all the linear
relations among Euler sums can be produced by FDS and EDS to be
defined below. In order to state it formally we need to adopt the
machinery in \cite[\S3]{IKZ}. We will use the same notations there
except that $\mathfrak H$ is replaced by $\fA$ and $y$ by $b$.
Then let $R$ be a commutative $\Q$-algebra with 1 and $Z_R$ is any
map from $\fA^0$ to $R$ such that the ``finite double shuffle''
(FDS) property holds:
  $$Z_R(w_1\sha w_2)=Z_R(w_1* w_2)=Z_R(w_1)Z_R(w_2).$$
We then extend $Z_R$ to $Z_R^\sha$ and $Z_R^*$ as before. Define
an $R$-module $R$-linear automorphism $\rho_R$ of $R[T]$ by
$$\rho_R(e^{Tu})=A_R(u)e^{Tu}$$
where
$$A_R(u)=\exp\left(\sum_{n=2}^\infty \frac{(-1)^n}{n}
 Z_R(a^{n-1} b)u^n\right) \in R[\![u]\!].$$
Similar to the situation for MZVs, we may define the
$\fA^0$-algebra isomorphisms
$$\reg_\sha^T:\fA_\sha^1 =\fA_\sha^0[b]\lra\fA_\sha^0[T],\qquad
 \reg_*^T:\fA_*^1 =\fA_*^0[b]\lra\fA_*^0[T],$$
which send $b$ to $T$. Composing these with the evaluation map
$T=0$ we get the maps $\reg_\sha$ and $\reg_*$.
\begin{conj} Let $(R,Z_R)$ be as above with the FDS property. Then
the following are equivalent:
\begin{itemize}
    \item[\emph{(i)}] $(Z_R^\sha-\rho_R\circ Z_R^*)(w)=0$ for all $w\in \fA^1$.
    \item[\emph{(ii)}]  $(Z_R^\sha-\rho_R\circ Z_R^*)(w)|_{T=0}=0$ for all $w\in \fA^1$.
   \item[\emph{(iii)}]  $ Z_R^\sha(w_1\sha w_0-w_1*w_2)=0$ for all $w_1\in
   \fA^1$ and all $w_0\in \fA^0$.
   \item[\emph{(iii$'$)}]  $ Z_R^*(w_1\sha w_0-w_1*w_2)=0$ for all $w_1\in
   \fA^1$ and all $w_0\in \fA^0$.
   \item[\emph{(iv)}]  $ Z(\reg_\sha(w_1\sha w_0-w_1*w_2))=0$ for all $w_1\in
   \fA^1$ and all $w_0\in \fA^0$.
   \item[\emph{(iv$'$)}]  $ Z(\reg_*(w_1\sha w_0-w_1*w_2))=0$ for all $w_1\in
   \fA^1$ and all $w_0\in \fA^0$.
   \item[\emph{(v)}]  $ Z(\reg_\sha(b^m*w))=0$ for all $m\ge 1$ and all $w\in \fA^0$.
   \item[\emph{(v$'$)}]  $ Z(\reg_*(b^m\sha w))=0$ for all $m\ge 1$ and all $w\in \fA^0$.
\end{itemize}
\end{conj}

If a map $Z_R:\fA^0\lra R$ satisfies the FDS and any one of the
equivalent conditions in the above conjecture then we say that
$Z_R$ have the extended double shuffle (EDS) property. Let
$R_{EDS}$ be the universal algebra (together with a map
$Z_{EDS}:\fA^0\lra R_{EDS}$) such that for every $\Q$-algebra $R$
and a map $Z_R:\fA^0\lra R$ satisfying EDS there always exists a
map $\varphi_R$ to make the following diagram commutative:
\begin{equation*}
\text{$\diagramcompileto{diagwt2}
\fA^0 \drto_{Z_R}\rto^-{Z_{EDS}}& R_{EDS} \dto^-{\varphi_R}\\
\ & R
\enddiagram$}
\end{equation*}
\begin{mainconj} 
The map $\phi_\R$ is injective, namely, the algebra of Euler sums
is isomorphic to $R_{EDS}.$
\end{mainconj}

If an Euler sum can be expressed by linear combination of the
products of Euler sums with lower weights then the Euler sum is
called \emph{reduced}. Broadhurst \cite{Br2} gives a conjecture on
the number of Euler sums in a minimal $\Q$-basis for reducing all
Euler sums to basic Euler sums. When considering only the linear
independence of Euler sums Broadhurst conjectures that the
$\Q$-dimension of weight $n$ Euler sum sums is given by the
Fibonacci numbers: $d_1=2,$ $d_3=3,$ $d_4=5$, $d_5=8$, and so on.
Zlobin \cite{Zl} further proposes the following precise version of
this conjecture.
\begin{conj} Every weight $n$ Euler sum is a $\Q$-linear
combination of the following Euler sums: $\zeta(\bar
b_1,b_2,\dots,b_r)$, where $b_j\in\{1,2\}$ and $\sum_{j=0}^r
b_j=n$.
\end{conj}

However, further computation suggests there may exist even subtler
structures. So we propose
\begin{mainconj} Let $n$ be a positive integer. Then there are
$\Q$-linearly independent Euler sums of weight $n$ such that every
Euler sum of weight $n$ is a $\Z$-linear combination of these
sums.
\end{mainconj}
We will denote $EZ_n$ (for ``Euler sums relations over $\Z$'') the
number of independent Euler sums of weight $n$ in the conjecture.
It is likely that $EZ_2=2,$ $EZ_3=3$, $EZ_4=5$ and $EZ_5=8$ which
are suggested by the computations in the next section, which agree
with Broadhurst's conjecture. In another paper \cite{Zh} we
investigate the relations between special values of multiple
polylogarithms at $m$th roots of unity for general $m$ and propose
a similar problem to Main Conjecture 2.

\section{The structure of Euler sums and some numerical
evidence}\label{sec:comp} We shall now use both FDS and EDS to
compute the relations between Euler sums of weight $<6$. Most of
the computations in this section are carried out by Maple. We have
checked the consistency of these relations with the many known
ones and verified numerically all the identities in the paper by
EZ-face \cite{EZface} with error smaller than $10^{-50}$. From
these numerical results we derived our Main Conjecture 2.

\begin{prop} All the weight two Euler sums can be expressed as
$\Z$-linear combinations of $\zeta(\bar2)$ and $\zeta(\bar1,1):$
\begin{equation*}
  \zeta(2)=-2\zeta(\bar2),\quad
 \zeta(\bar1, \bar1)=\zeta(\bar2)+\zeta(\bar1,1)
\end{equation*}\end{prop}
\begin{proof} It is easy to see from EDS that
\begin{equation*}
\zeta(2)=-2\zeta(\bar1, \bar1) + 2 \zeta(\bar1, 1),\qquad
 \zeta(\bar2)=-\zeta(\bar1,1) + \zeta(\bar1, \bar1).
\end{equation*}
\end{proof}
\begin{rem} From the proposition and a stuffle relation we get
 $$2\zeta(\bar1,1)=2\zeta(\bar1,
 \bar1)-2\zeta(\bar2)=\zeta(\bar1)^2=\ln(2)^2.$$
Hence it is apparent that $\zeta(2)$ and $\zeta(\bar1,1)$ are
linearly independent over $\Q$ which verifies the Main
Conjecture~1 in this case.
\end{rem}

\begin{prop}\label{prop:wt3EDS}
We can express all weight three Euler sums as $\Z$-linear
combinations of $\zeta(\bar2,1),$  $\zeta(\bar1,1,1)$ and
$\zeta(\bar1,2)$:
 {\setlength\arraycolsep{1pt}
$$
\ \hskip-.2cm
\begin{array}{rrcrlrl}
\zeta(3) =& 8\zeta(\bar2,1)&,\\
\zeta(\bar3) =& -6\zeta(\bar2,1)&,\\
\zeta(2,1) =& 8\zeta(\bar2,1)&,\\
\zeta(2,\bar1) =& 2\zeta(\bar2,1)&-&3\zeta(\bar1,2)&,\\
\zeta(\bar2,\bar1) =&3\zeta(\bar1,2)&-&7\zeta(\bar2,1)&,\\
\zeta(\bar1,\bar2) =& -2\zeta(\bar1,2)&+&\zeta(\bar2,1)&,\\
\zeta(\bar1,1,\bar1) =& \zeta(\bar2,1)&+&\zeta(\bar1,1,1)&,\\
\zeta(\bar1,\bar1,1) =&\zeta(\bar1,2)&-&5\zeta(\bar2,1)&+&\zeta(\bar1,1,1)&,\\
\zeta(\bar1,\bar1,\bar1) =& \zeta(\bar1,2)&+&\zeta(\bar1,1,1)&.
\end{array}$$}
\end{prop}
\begin{proof}
When weight is three, by only DS we have
\begin{alignat*}{3}
\ &\zeta(\bar1,1,\bar1)+2\zeta(\bar1,\bar1,1)+\zeta(\bar1,\bar2)+\zeta(2,1)-3\zeta(\bar1,1,1)&=0,\\
\ &2\zeta(\bar1,\bar1,\bar1)+\zeta(\bar1,2)+\zeta(2,\bar1)-2\zeta(\bar1,1,\bar1)&=0,\\
\ &\zeta(\bar2,\bar1)+\zeta(\bar1,\bar2)+\zeta(3)-2\zeta(\bar2,1)-\zeta(\bar1,2)&=0,\\
\
&\zeta(\bar1,2)+\zeta(\bar3)-\zeta(\bar2,\bar1)-\zeta(\bar1,\bar2)&=0.
\end{alignat*}
These are far from enough to prove the proposition. But by EDS we
have five more relations:
\begin{alignat*}{3}
\
&\zeta(\bar3)+2\zeta(\bar2,1)+\zeta(\bar1,2)+2\zeta(\bar1,1,1)-\zeta(2,\bar1)+\zeta(\bar1)\zeta(2)-2\zeta(\bar1,\bar1,1)&=0,\\
\ &\zeta(\bar1,1,\bar1)-\zeta(\bar2,1)-\zeta(\bar1,2)-2\zeta(\bar1,1,1)+\zeta(\bar1,\bar1,\bar1)&=0,\\
\ &\zeta(\bar1,\bar1,1)-\zeta(\bar2,\bar1)-\zeta(\bar1,\bar2)-\zeta(\bar1,1,\bar1)&=0,\\
\ &\zeta(\bar2,\bar1)-\zeta(\bar3)-\zeta(\bar2,1)+\zeta(2,\bar1)&=0,\\
\ &\zeta(2,1)-\zeta(3)&=0.
\end{alignat*}
Now the proposition follows from the stuffle relation:
$\zeta(\bar1)\zeta(2)=\zeta(\bar3)+\zeta(2,\bar1)+\zeta(\bar1,2).$
\end{proof}

\begin{rem} By our Main Conjecture 1 there should be no further
linear relations among $\zeta(\bar2,1),$ $\zeta(\bar1,1,1)$ and
$\zeta(\bar1,2)$ which gives $EZ_3=3$. This is easy to see to be
equivalence to the linear independence of $\zeta(3),$
$\zeta(\bar1)\zeta(2)$ and $\zeta(\bar1,\bar1,1)$.
\end{rem}

The previous two propositions and the following two results show
that if weight $<6$ then both Broadhurst-Zlobin Conjecture and our
Main Conjecture 2 are true.
\begin{prop}\label{prop:wt4EDS}
All weight four Euler sums are $\Z$-linear combinations of
$A=\zeta(\bar2,1,1)$, $B=\zeta(\bar2,2),$ $C=\zeta(\bar1,2,1),$
$D=\zeta(\bar1,1,2),$ and $E=\zeta(\bar1,1,1,1)$. For length one
and two:
 {\setlength\arraycolsep{1pt}
$$
\ \hskip-.2cm
\begin{array}{rrlrlrlrl}
\zeta(4)=&64A&+&16B&, \\ \zeta(\bar4)=&-56A&-&14B&, \\
\zeta(3,1)=&16A&+&4B&, \\ \zeta(3,\bar1)=&118A&+&19B&+&14C&, \\
\zeta(2,2)=&48A&+&12B&, \\ \zeta(\bar3,1)=&10A&+&2B&, \\
\zeta(\bar3,\bar1)=&-140A&-&24B&-&14C&, \\
\zeta(2,\bar2)=&-24A&-&7B&, \\ \zeta(\bar2,\bar2)=&-12A&-&3B&, \\
\zeta(\bar1,3)=&-38A&-&5B&-&6C&, \\
\zeta(\bar1,\bar3)=&58A&+&8B&+&8C&.
\\ \end{array}$$}
For length three:
 {\setlength\arraycolsep{1pt}
$$
\ \hskip-.2cm
\begin{array}{rrlrlrlrl}
\zeta(2,1,1)=&64A&+&16B&, \\
\zeta(2,1,\bar1)=&16A&+&2B&+&6C&+&3D&, \\
\zeta(2,\bar1,1)=&22A&+&3B&+&C&-&3D&, \\
\zeta(2,\bar1,\bar1)=&100A&+&13B&+&9C&-&6D&, \\
\zeta(\bar2,1,\bar1)=&91A&+&14B&+&8C&-&3D&, \\
\zeta(\bar2,\bar1,1)=&-161A&-&26B&-&15C&+&3D&, \\
\zeta(\bar2,\bar1,\bar1)=&-101A&-&14B&-&9C&+&6D&, \\
\zeta(\bar1,2,\bar1)=&-102A&-&14B&-&8C&+&6D&, \\
\zeta(\bar1,\bar2,1)=&69A&+&11B&+&8C&, \\
\zeta(\bar1,\bar2,\bar1)=&63A&+&8B&+&3C&-&6D&, \\
\zeta(\bar1,1,\bar2)=&21A&+&3B&+&C&-&2D&, \\
\zeta(\bar1,\bar1,\bar2)=&A&+&2B& && +&D&.
\end{array}$$}
For length four,
 {\setlength\arraycolsep{1pt}
$$
\ \hskip-.2cm
\begin{array}{rrlrlrlrlrlrl}
\zeta(\bar1,1,1,\bar1)=&A&&&&&&&+&E&, \\
\zeta(\bar1,1,\bar1,1)=&11A&+&2B&+&C&&&+&E&, \\
\zeta(\bar1,1,\bar1,\bar1)=&&&&&C&&&+&E&, \\
\zeta(\bar1,\bar1,1,1)=&-83A&-&16B&-&5C&+&D&+&E&, \\
\zeta(\bar1,\bar1,1,\bar1)=&-38A&-&5B&-&5C&+&D&+&E&, \\
\zeta(\bar1,\bar1,\bar1,1)=&&&&&&&D&+&E&, \\
\zeta(\bar1,\bar1,\bar1,\bar1)=&A&+&B&&&+&D&+&E&.
\end{array}$$}
\end{prop}

The next proposition shows that the $\Q$-basis conjectured by
Zlobin can not be chosen as the $\Z$-linear basis in general.
\begin{prop}\label{prop:wt5EDS}
All weight five Euler sums are $\Q$-linear combinations of
$\zeta(\bar1,1,1,1,1)$, $\zeta(\bar1,1,2,1)$,
$\zeta(\bar2,1,1,1)$, $\zeta(\bar1,1,1,2)$, $\zeta(\bar1,2,1,1)$,
$\zeta(\bar2,1,2)$, $\zeta(\bar2,2,1)$ and $\zeta(\bar1,2,2)$. For
example
$$\zeta(3,1,1)=-\frac{448}{39}\zeta(\bar2,1,1,1)-\frac{112}{39}\zeta(\bar2,2,1)-\frac{48}{13}\zeta(\bar2,1,2) .$$
Furthermore, all weight five Euler sums are $\Z$-linear
combinations of {\setlength\arraycolsep{1pt}
$$
\ \hskip-.2cm
\begin{array}{rlrlrlrl}
A=&\zeta(\bar1,\bar1,\bar1,2),& B=&\zeta(\bar2,1,\bar1,\bar1),&
C=&\zeta(\bar1,1,\bar1,\bar2),& D=&\zeta(\bar2,1,1,1),\\
E=&\zeta(\bar1,\bar1,\bar1,1,1),& F=&\zeta(2,2,\bar1),&
G=&\zeta(\bar1,1,\bar1,1,\bar1),&
H=&\zeta(\bar1,\bar1,\bar1,\bar1,\bar1).
\end{array}$$}
Eor length one and two: {\setlength\arraycolsep{1pt}
$$
\ \hskip-.2cm
\begin{array}{rrlrlrlrlrlrlrlrlrl}
\zeta(5)=&-13504A&+&1856B&-&1344C&+&26880D&-&18752E&-&640F&-&31552G&+&50304H&,
\\
\zeta(\bar5)=&12660A&-&1740B&+&1260C&-&25200D&+&17580E&+&600F&+&29580G&-&47160H&,
\\
\zeta(4,1)=&-9808A&+&1344B&-&944C&+&19632D&-&13648E&-&464F&-&22848G&+&36496H&,
\\
\zeta(4,\bar1)=&-14918A&+&2044B&-&1434C&+&29862D&-&20758E&-&704F&-&34748G&+&55506H&,
\\
\zeta(\bar4,1)=&3638A&-&498B&+&346C&-&7296D&+&5066E&+&172F&+&8466G&-&13532H&,
\\
\zeta(\bar4,\bar1)=&19862A&-&2722B&+&1914C&-&39744D&+&27634E&+&938F&+&46274G&-&73908H&,
\\
\zeta(3,2)=&22672A&-&3104B&+&2160C&-&45456D&+&31568E&+&1072F&+&52768G&-&84336H&,
\\
\zeta(3,\bar2)=&4562A&-&626B&+&446C&-&9108D&+&6342E&+&216F&+&10642G&-&16984H&,
\\
\zeta(\bar3,2)=&-6552A&+&898B&-&632C&+&13110D&-&9116E&-&310F&-&15266G&+&24382H&,
\\
\zeta(\bar3,\bar2)=&-17848A&+&2444B&-&1704C&+&35772D&-&24848E&-&844F&-&41548G&+&66396H&,
\\
\zeta(2,3)=&-26368A&+&3616B&-&2560C&+&52704D&-&36672E&-&1248F&-&61472G&+&98144H&,
\\
\zeta(2,\bar3)=&6792A&-&934B&+&680C&-&13506D&+&9428E&+&322F&+&15878G&-&25306H&,
\\
\zeta(\bar2,\bar3)=&24902A&-&3412B&+&2394C&-&49854D&+&34654E&+&1178F&+&58004G&-&92658H&,
\\
\zeta(\bar2,3)=&-8622A&+&1182B&-&834C&+&17244D&-&11994E&-&408F&-&20094G&+&32088H&,
\\
\zeta(\bar1,4)=&5266A&-&720B&+&494C&-&10582D&+&7338E&+&248F&+&12240G&-&19578H&,
\\
\zeta(\bar1,\bar4)=&-8990A&+&1230B&-&850C&+&18044D&-&12522E&-&424F&-&20910G&+&33432H&.
\end{array}$$}
For length three, {\setlength\arraycolsep{1pt}
$$
\ \hskip-.2cm
\begin{array}{rrlrlrlrlrlrlrlrlrl}
\zeta(3,1,1)=&-9808A&+&1344B&-&944C&+&19632D&-&13648E&-&464F&-&22848G&+&36496H&,
\\
\zeta(3,1,\bar1)=&-5314A&+&725B&-&500C&+&10677D&-&7402E&-&250F&-&12339G&+&19741H&,
\\
\zeta(3,\bar1,1)=&-2257A&+&312B&-&225C&+&4489D&-&3137E&-&108F&-&5290G&+&8427H&,
\\
\zeta(3,\bar1,\bar1)=&-7299A&+&1005B&-&713C&+&14566D&-&10151E&-&347F&-&17057G&+&27208H&,
\\
\zeta(\bar3,1,1)=&4482A&-&614B&+&430C&-&8974D&+&6238E&+&212F&+&10438G&-&16676H&,
\\
\zeta(\bar3,1,\bar1)=&9570A&-&1308B&+&908C&-&19204D&+&13328E&+&452F&+&22250G&-&35578H&,
\\
\zeta(\bar3,\bar1,1)=&12462A&-&1710B&+&1204C&-&24924D&+&17338E&+&590F&+&29056G&-&46394H&,
\\
\zeta(\bar3,\bar1,\bar1)=&-4288A&+&582B&-&396C&+&8646D&-&5978E&-&201F&-&9922G&+&15900H&,
\\
\zeta(2,2,1)=&22672A&-&3104B&+&2160C&-&45456D&+&31568E&+&1072F&+&52768G&-&84336H&,
\\
\zeta(2,\bar2,1)=&3025A&-&414B&+&287C&-&6065D&+&4213E&+&143F&+&7038G&-&11251H&,
\\
\zeta(2,\bar2,\bar1)=&6421A&-&881B&+&627C&-&12818D&+&8927E&+&303F&+&14977G&-&23904H&,
\\
\zeta(2,1,2)=&-26368A&+&3616B&-&2560C&+&52704D&-&36672E&-&1248F&-&61472G&+&98144H&,
\\
\zeta(2,\bar1,2)=&2206A&-&302B&+&210C&-&4428D&+&3074E&+&104F&+&5134G&-&8208H&,
\\
\zeta(2,1,\bar2)=&7958A&-&1093B&+&786C&-&15861D&+&11056E&+&377F&+&18581G&-&29637H&,
\\
\zeta(2,\bar1,\bar2)=&23513A&-&3221B&+&2255C&-&47094D&+&32727E&+&1113F&+&54757G&-&87484H&,
\\
\zeta(\bar2,2,1)=&-12813A&+&1755B&-&1227C&+&25664D&-&17835E&-&606F&-&29835G&+&47670H&,
\\
\zeta(\bar2,2,\bar1)=&-20468A&+&2804B&-&1964C&+&40988D&-&28488E&-&968F&-&47668G&+&76156H&,
\\
\zeta(\bar2,\bar2,1)=&-5477A&+&750B&-&523C&+&10977D&-&7625E&-&259F&-&12750G&+&20375H&,
\\
\zeta(\bar2,\bar2,\bar1)=&12308A&-&1686B&+&1180C&-&24654D&+&17132E&+&582F&+&28662G&-&45794H&,
\\
\zeta(\bar2,1,2)=&12622A&-&1729B&+&1210C&-&25281D&+&17568E&+&597F&+&29393G&-&46961H&,
\\
\zeta(\bar2,1,\bar2)=&-3065A&+&420B&-&295C&+&6135D&-&4265E&-&145F&-&7140G&+&11405H&,
\\
\zeta(\bar2,\bar1,2)=&-14047A&+&1923B&-&1337C&+&28170D&-&19561E&-&665F&-&32691G&+&52252H&,
\\
\zeta(\bar2,\bar1,\bar2)=&-9411A&+&1290B&-&909C&+&18831D&-&13095E&-&445F&-&21930G&+&35025H&,
\\
\zeta(\bar1,3,1)=&123A&-&17B&+&13C&-&242D&+&171E&+&6F&+&289G&-&460H&,
\\
\zeta(\bar1,3,\bar1)=&-11820A&+&1614B&-&1120C&+&23726D&-&16460E&-&557F&-&27466G&+&43926H&,
\\
\zeta(\bar1,\bar3,1)=&6380A&-&874B&+&612C&-&12776D&+&8880E&+&302F&+&14858G&-&23738H&,
\\
\zeta(\bar1,\bar3,\bar1)=&12610A&-&1722B&+&1194C&-&25312D&+&17560E&+&594F&+&29302G&-&46862H&,
\\
\zeta(\bar1,2,2)=&-190A&+&26B&-&18C&+&384D&-&266E&-&9F&-&442G&+&708H&,
\\
\zeta(\bar1,2,\bar2)=&13726A&-&1880B&+&1314C&-&27494D&+&19106E&+&649F&+&31960G&-&51066H&,
\\
\zeta(\bar1,\bar2,2)=&-13631A&+&1867B&-&1305C&+&27302D&-&18973E&-&644F&-&31739G&+&50712H&,
\\
\zeta(\bar1,\bar2,\bar2)=&599A&-&82B&+&57C&-&1203D&+&835E&+&28F&+&1394G&-&2229H&,
\\
\zeta(\bar1,1,3)=&-3186A&+&435B&-&300C&+&6399D&-&4438E&-&150F&-&7401G&+&11839H&,
\\
\zeta(\bar1,1,\bar3)=&-2732A&+&376B&-&268C&+&5448D&-&3798E&-&130F&-&6384G&+&10182H&,
\\
\zeta(\bar1,\bar1,3)=&20431A&-&2799B&+&1969C&-&40888D&+&28427E&+&966F&+&47591G&-&76018H&,
\\
\zeta(\bar1,\bar1,\bar3)=&-7808A&+&1070B&-&758C&+&15608D&-&10858E&-&369F&-&18196G&+&29054H&.
\end{array}$$}
For length four, {\setlength\arraycolsep{1pt}
$$
\ \hskip-.2cm
\begin{array}{rrlrlrlrlrlrlrlrlrl}
\zeta(2,1,1,1)=&-13504A&+&1856B&-&1344C&+&26880D&-&18752E&-&640F&-&31552G&+&50304H&,
\\
\zeta(2,1,1,\bar1)=&-11109A&+&1518B&-&1044C&+&22320D&-&15477E&-&523F&-&25812G&+&41289H&,
\\
\zeta(2,1,\bar1,1)=&1174A&-&158B&+&101C&-&2395D&+&1642E&+&54F&+&2691G&-&4333H&,
\\
\zeta(2,1,\bar1,\bar1)=&14927A&-&2044B&+&1431C&-&29899D&+&20773E&+&705F&+&34745G&-&55518H&,
\\
\zeta(2,\bar1,1,1)=&-2712A&+&371B&-&258C&+&5439D&-&3777E&-&128F&-&6306G&+&10083H&,
\\
\zeta(2,\bar1,1,\bar1)=&-14828A&+&2030B&-&1419C&+&29709D&-&20641E&-&701F&-&34517G&+&55158H&,
\\
\zeta(2,\bar1,\bar1,1)=&-7120A&+&977B&-&681C&+&14262D&-&9911E&-&337F&-&16585G&+&26496H&,
\\
\zeta(2,\bar1,\bar1,\bar1)=&11204A&-&1534B&+&1074C&-&22440D&+&15595E&+&530F&+&26096G&-&41691H&,
\\
\zeta(\bar2,\bar1,1,1)=&8717A&-&1197B&+&847C&-&17415D&+&12122E&+&412F&+&20334G&-&32456H&,
\\
\zeta(\bar2,\bar1,1,\bar1)=&-8511A&+&1162B&-&806C&+&17085D&-&11852E&-&401F&-&19775G&+&31627H&,
\\
\zeta(\bar2,\bar1,\bar1,1)=&3432A&-&470B&+&327C&-&6882D&+&4779E&+&162F&+&7980G&-&12759H&,
\\
\zeta(\bar2,\bar1,\bar1,\bar1)=&-652A&+&89B&-&66C&+&1296D&-&905E&-&31F&-&1531G&+&2436H&,
\\
\zeta(\bar2,1,\bar1,1)=&8659A&-&1183B&+&822C&-&17376D&+&12059E&+&409F&+&20134G&-&32193H&,
\\
\zeta(\bar2,1,1,\bar1)=&3571A&-&490B&+&344C&-&7145D&+&4969E&+&169F&+&8322G&-&13291H&,
\\
\zeta(\bar1,2,1,1)=&190A&-&26B&+&18C&-&384D&+&265E&+&9F&+&442G&-&707H&,
\\
\zeta(\bar1,2,1,\bar1)=&190A&-&25B&+&18C&-&385D&+&265E&+&9F&+&442G&-&707H&,
\\
\zeta(\bar1,2,\bar1,1)=&-27776A&+&3801B&-&2654C&+&55667D&-&38668E&-&1313F&-&64648G&+&103316H&,
\\
\zeta(\bar1,2,\bar1,\bar1)=&-19006A&+&2604B&-&1828C&+&38048D&-&26452E&-&900F&-&44288G&+&70740H&,
\\
\zeta(\bar1,\bar2,1,1)=&-2407A&+&330B&-&233C&+&4812D&-&3347E&-&113F&-&5610G&+&8957H&,
\\
\zeta(\bar1,\bar2,1,\bar1)=&-202A&+&32B&-&34C&+&353D&-&274E&-&12F&-&533G&+&807H&,
\\
\zeta(\bar1,\bar2,\bar1,1)=&6507A&-&893B&+&631C&-&13009D&+&9054E&+&309F&+&15184G&-&24238H&,
\\
\zeta(\bar1,\bar2,\bar1,\bar1)=&31628A&-&4333B&+&3038C&-&63328D&+&44021E&+&1497F&+&73681G&-&117702H&,
\\
\zeta(\bar1,1,2,1)=&-122A&+&17B&-&13C&+&242D&-&170E&-&6F&-&288G&+&458H&,
\\
\zeta(\bar1,1,2,\bar1)=&3310A&-&453B&+&314C&-&6641D&+&4609E&+&156F&+&7692G&-&12301H&,
\\
\zeta(\bar1,1,\bar2,\bar1)=&6195A&-&850B&+&600C&-&12383D&+&8619E&+&294F&+&14454G&-&23073H&,
\\
\zeta(\bar1,1,\bar2,1)=&2888A&-&394B&+&272C&-&5803D&+&4023E&+&136F&+&6706G&-&10729H&,
\\
\zeta(\bar1,\bar1,2,1)=&-7433A&+&1019B&-&711C&+&14888D&-&10348E&-&352F&-&17315G&+&27663H&,
\\
\zeta(\bar1,\bar1,2,\bar1)=&6793A&-&932B&+&658C&-&13586D&+&9453E&+&322F&+&15848G&-&25301H&,
\\
\zeta(\bar1,\bar1,\bar2,1)=&-313A&+&43B&-&30C&+&626D&-&437E&-&15F&-&730G&+&1167H&,
\\
\zeta(\bar1,\bar1,\bar2,\bar1)=&-18914A&+&2592B&-&1822C&+&37855D&-&26321E&-&895F&-&44073G&+&70394H&,
\\
\zeta(\bar1,1,1,2)=&191A&-&26B&+&18C&-&384D&+&267E&+&9F&+&442G&-&709H&,
\\
\zeta(\bar1,1,1,\bar2)=&-2521A&+&345B&-&240C&+&5054D&-&3510E&-&119F&-&5864G&+&9374H&,
\\
\zeta(\bar1,1,\bar1,2)=&13126A&-&1798B&+&1257C&-&26291D&+&18271E&+&621F&+&30567G&-&48838H&,
\\
\zeta(\bar1,\bar1,1,2)=&13312A&-&1826B&+&1295C&-&26595D&+&18512E&+&630F&+&31043G&-&49555H&,
\\
\zeta(\bar1,\bar1,1,\bar2)=&-4812A&+&661B&-&475C&+&9593D&-&6687E&-&228F&-&11237G&+&17924H&,
\\
\zeta(\bar1,\bar1,\bar1,\bar2)=&-13127A&+&1798B&-&1258C&+&26291D&-&18271E&-&621F&-&30565G&+&48836H&.
\end{array}$$}
For length five, {\setlength\arraycolsep{1pt}
$$
\ \hskip-.2cm
\begin{array}{rrlrlrlrlrlrlrlrlrl}
\zeta(\bar1,1,1,1,1)=&-191A&+&26B&-&18C&-&442G&+&384D&-&266E&-&9F&+&709H&,
\\
\zeta(\bar1,1,1,1,\bar1)=&-191A&+&26B&-&18C&+&385D&-&266E&-&9F&-&442G&+&709H&,
\\
\zeta(\bar1,1,1,\bar1,1)=&4481A&-&614B&+&430C&-&8973D&+&6237E&+&212F&+&10438G&-&16674H&,
\\  \zeta(\bar1,1,1,\bar1,\bar1)=&-A&&&&&&&-&E&&&&&+&2H&, \\
\zeta(\bar1,1,\bar1,1,1)=&-4693A&+&643B&-&451C&+&9395D&-&6531E&-&222F&-&10930G&+&17462H&,
\\
\zeta(\bar1,1,\bar1,\bar1,1)=&-313A&+&43B&-&31C&-&730G&+&626D&-&436E&-&15F&+&1167H&,
\\
\zeta(\bar1,1,\bar1,\bar1,\bar1)=&-13126A&+&1798B&-&1258C&+&26291D&-&18271E&-&621F&-&30565G&+&48837H&,
\\
\zeta(\bar1,\bar1,1,1,1)=&7496A&-&1031B&+&747C&-&14915D&+&10408E&+&355F&+&17522G&-&27929H&,
\\
\zeta(\bar1,\bar1,1,1,\bar1)=&2081A&-&285B&+&194C&-&4183D&+&2901E&+&98F&+&4840G&-&7740H&,
\\
\zeta(\bar1,\bar1,1,\bar1,1)=&-3308A&+&452B&-&313C&+&6641D&-&4607E&-&156F&-&7689G&+&12297H&,
\\
\zeta(\bar1,\bar1,1,\bar1,\bar1)=&-12121A&+&1660B&-&1164C&+&24269D&-&16868E&-&573F&-&28225G&+&45094H&,
\\
\zeta(\bar1,\bar1,\bar1,1,\bar1)=&12622A&-&1729B&+&1210C&-&25280D&+&17569E&+&597F&+&29393G&-&46961H&,
\\
\zeta(\bar1,\bar1,\bar1,\bar1,1)=&10552A&-&1445B&+&1008C&-&21144D&+&14690E&+&499F&+&24565G&-&39254H&.
\end{array}$$}
\end{prop}

\section{A family of Euler sum identities}
In this section we shall prove the following
\begin{thm}  \label{thm:iden} {\em For every positive integer} $n$
$$\zeta(\{3\}^n)=8^n\zeta(\{\ol2,1\}^n).$$
\end{thm}
First we can rephrase our identities using words in $\fA^0$, which
was pointed out to us by D.~Bradley. For any positive integer $i$
define the $i$-th cut of a word $l_1\dots l_m$ ($l_i$ are letters)
to be a pair of words given by
\begin{equation*}
  \text{Cut}_i[l_1l_2,\dots l_m]=
  \begin{cases}
  \big[(l_1,l_2,\dots,l_i),(l_{i+1},\dots,l_m)\big] &\text{if $i$ is
  odd,}\\
   \big[(l_i,\dots,l_2,l_1),(l_{i+1},\dots,l_m)\big] &\text{if $i$ is
   even},
  \end{cases}
\end{equation*}
for $i=0,\cdots, m$. Here by convention for empty word $\bf 1$ we
have $[w,{\bf 1}]=[{\bf 1},w]=w$. For any two words $w_1,w_2$, set
 $$\sha[w_1,w_2]= w_1\sha w_2,\quad\text{and} \quad *[w_1,w_2]= w_1* w_2.$$
Then we can define the composites $\sha_i=\sha\circ \text{Cut}_i$,
$*_i=*\circ \text{Cut}_i$, and the difference $\gD_i=\sha_i-*_i$.
\begin{thm} \label{thm:key}
For a positive real number $x$ let $[x]$ and $\{x\}$ be the
integral part and the fractional part of $x$, respectively. For
any two words $l_1$ and $l_2$ define the $\star$-concatenation by
setting $l_1\star l_2=l_1 l_2$ except that
 $$b\star b=bc,\quad\text{and}\quad c\star c=cb.$$
Then for every positive integer $n$ the following holds in
$\fA^0$:
 \begin{equation}\label{equ:key}
2^n(ac^2ab^2)^{[n/2]}(ac^2)^{2\{n/2\}}=(a^2(b+c))^n
+\sum_{i=0}^{2n} (-1)^{n-i}\gD_i\big((cd)^{\star n}).
\end{equation}
Here $d=a(b+c)$ is regarded as one letter when we do the cuts
first, retaining the $\star$-concatenation.
\end{thm}
\begin{rem} (1) Observe that $(b+c)\star c=(b+c)\star b=bc+cb$ and
therefore
$$Z((cd)^{\star n})=\sum_{t_1,\dots,t_n\in\{2,\bar
 2\}}\zeta\big(\bar1,t_1,\bar1,t_2,\dots,\bar1,t_n\big).$$
The $\star$-concatenation appears because neither $c^2$ nor $b^2$
can appear in any of the Euler sums in the above sum. We should
keep this in mind because the operator Cut$_i$ will lead to some
order reversals which also should obey this condition.

(2) As pointed out by D.~Bradley the theorem is quite similar to
\cite[Lemma 3.1]{BoB} in spirit although they are not the same. Is
there any relation between them?
\end{rem}
It is easy to verify that for any positive integer $n$
$$\zeta(\{\bar2,1\}^n)=(ac^2ab^2)^{[n/2]}(ac^2)^{2\{n/2\}}.$$
and
$$\zeta(\{3\}^n)=Z\big((a^2b)^n\big).$$
On the other hand an integral substitution $t\to t^2$ yields (see
\cite[(5.14)]{BdB})
 \begin{equation}\label{subs}
 \zeta(\{3\}^n)=4^nZ\big((a^2(b+c))^n\big).
\end{equation}
This also follows quickly from the following special case of the
distribution relation of multiple polylogarithms (see
\cite[(2.5)]{Zh}):
  \begin{align*}
Z\big((a^2(b+c))^n\big)=& \sum_{1\le i_1<\cdots<i_j\le n}\
\sum_{k_1>\cdots>k_n>0}
 \frac{(-1)^{k_{i_1}+\cdots+k_{i_j}}}{(k_1\cdots k_n)^3}\\
 =&\sum_{k_1>\cdots>k_n>0}
 \frac{(1+(-1)^{k_1})\cdots(1+(-1)^{k_n})}{(k_1\cdots
 k_n)^3}=\frac{1}{4^n} \zeta(\{3\}^n).
\end{align*}
\begin{rem}
 From Maple computation we notice that \eqref{subs} can not be derived from FDS in
general. But we believe it is a consequence of some EDS from
Prop.~\ref{prop:eDS}. We plan to study this problem and EDS in
more details in the future.
\end{rem}

Now we can multiply $4^n$ on both sides of \eqref{equ:key} and
then apply $Z$. From Prop.~\ref{*homo} and Prop.~\ref{shahomo} we
see immediately that our Main Theorem follows.

To prove Theorem \ref{thm:key} we need two separate identities
involving stuffles and shuffles respectively.
\begin{prop}\label{prop:stuffle} For every positive integer $n$
 \begin{equation}\label{equ:stuffle}
 \sum_{i=0}^{2n} (-1)^i*_i\big((cd)^{\star n})=(-1)^n(a^2(b+c))^n.
\end{equation}
\end{prop}
\begin{proof}
We proceed by induction on $n$. When $n=1$ the left hand side of
\eqref{equ:stuffle} is
 \begin{align*}
 cd-c*d+d\star c
 =&\gam_1(\gam_2+\gb_2)-\gam_1*(\gam_2+\gb_2)+a(bc+cb) \\
 =&-(\gam_2+\gb_2)\mal\gam_1-\gb_3-\gam_3+\gam_2\gb_1+\gb_2\gam_1
 =-\gb_3-\gam_3.
\end{align*}
This is exactly the right hand side $-a^2(b+c)$. Now assume that
identity \eqref{equ:stuffle} holds up to $n-1$ for some $n\ge 2$.
Set $\gam=\gam_1=c$, $d=\gb_2+\gam_2$ and $d_3=\gb_3+\gam_3$. Then
$d\star \gam=\gb_2\gam+\gam_2\gb_1$. In the rest of the paper we
set $\mal=\mal_\gam$. Note that $\mal(d)=d$ and $\mal(d\star
\gam)=d\star \gam$. Hence by the recursive definition of the
stuffle product \eqref{equ:defnstuffle}
 \allowdisplaybreaks{
  \begin{align*}
\ &\sum_{i=0}^{2n} (-1)^i*_i\big((cd)^{\star n}) \\
= &\sum_{j=0}^{n}
 (d\star \gam)^j*(\gam d)^{\star (n-j)}-\sum_{j=1}^{n}
 \Big(\gam(d\star \gam)^{j-1}\Big)*\Big(d\star (\gam d)^{(n-j)}\Big)\\
 =&\sum_{j=0}^{n}
 \Big(\gb_2\gam(d\star \gam)^{j-1}\Big)*(\gam d)^{\star (n-j)}
 -\sum_{j=1}^{n}
 \Big(\gam(d\star \gam)^{j-1}\Big)*\Big(\gb_2\gam(d\star \gam )^{(n-j)}d\Big)\\
 +&\sum_{j=0}^{n}
 \Big(\gam_2\gb_1(d\star \gam)^{j-1}\Big)*(\gam d)^{\star (n-j)}
 -\sum_{j=1}^{n}
 \Big(\gam(d\star \gam)^{j-1}\Big)*\Big(\gam_2\gb_1 (d\star \gam )^{(n-j)}d\Big)\\
 =&\sum_{j=1}^{n-1}
 \left\{\gb_2\Big(\big(\gam(d\star \gam)^{j-1}\big)*(\gam d)^{\star (n-j)}\Big)
 +\gam\mal\Big(\big(\gb_2\gam(d\star \gam)^{j-1}\big)*\big((d\star \gam )^{(n-j-1)}d\big)\Big)\right.\\
 \ & \ \ \ +\left.\gam_3\mal\Big(\big(\gam(d\star \gam)^{j-1}\big)*\big((d\star \gam
 )^{(n-j-1)}d\big)\Big)\right \}+(\gam d)^{\star n}+(d\star \gam )^n\\
 -&\sum_{j=1}^{n}
 \left\{\gam\mal\Big((d\star \gam)^{j-1}*\big(\gb_2\gam(d\star \gam )^{(n-j)}d\big)\Big)
 +\gb_2\Big(\big(\gam(d\star \gam)^{j-1}\big)*\big(\gam(d\star \gam )^{(n-j)}d\big)\Big)\right.\\
\ & \ \ \ + \left.\gam_3\mal\Big((d\star
\gam)^{j-1}*\big(\gam(d\star \gam
  )^{(n-j)}d\big)\Big)\right \}\\
+&\sum_{j=1}^{n-1}
 \left\{\gam_2\mal\Big(\big(\gam(d\star \gam)^{j-1}\big)*(\gam d)^{\star (n-j)}\Big)
 +\gam\mal\Big(\big(\gam_2\gb_1(d\star \gam)^{j-1}\big)*
 \big((d\star \gam)^{(n-j-1)}d\big)\Big)\right.\\
\ & \ \ \ +\left.\gb_3\Big(\big(\gam(d\star \gam)^{j-1}\big)*
 \big((d\star \gam)^{(n-j-1)}d\big)\Big)\right\}\\
 -&\sum_{j=1}^{n}\left\{
 \gam\mal\Big((d\star \gam)^{j-1}*\big(\gam_2\gb_1 (d\star \gam )^{(n-j)}d\big)\Big)
 +  \gam_2\mal\Big(\big(\gam(d\star \gam)^{j-1}\big)*\big(\gam (d\star \gam )^{(n-j)}d\big)\Big)\right.\\
 \ & \ \ \ + \left.\gam_3\mal\Big((d\star \gam)^{j-1}*
 \big(\gam (d\star \gam )^{(n-j)}d\big)\Big)\right\}.
\end{align*}}
Converting $\gb_2\gam+\gam_2\gb_1$ back to $d\star \gam$ and
cancelling all the terms without $\gam_3$ or $\gb_3$ we get
$$\sum_{i=0}^{2n} (-1)^i*_i\big((cd)^{\star n})
= d_3\left\{ \sum_{j=1}^{n-1} \big(\gam(d\star
\gam)^{j-1}\big)*\big( (d\star \gam )^{(n-j-1)}d\big)
-\sum_{j=1}^{n} (d\star \gam)^{j-1}*
 \big(\gam (d\star \gam )^{(n-j)}d\big)\right\}$$
By induction assumption the expression in the last big curly
bracket is $(-1)^n(a^2(b+c))^{n-1}.$ This proves the proposition
since $d_3=a^2(b+c)$.
\end{proof}

\begin{prop} \label{prop:shuffle} For every positive integer $n$
 \begin{equation}\label{equ:shuffle}
\sum_{i=0}^{2n} (-1)^i\sha_i\big((cd)^{\star n})=(-2)^n
 (ac^2ab^2)^{[n/2]}(ac^2)^{2\{n/2\}}
\end{equation}
and
 \begin{equation}\label{equ:shuffleb}
\sum_{i=0}^{2n} (-1)^i\sha_i\big((bd)^{\star n})=(-2)^n
 (ab^2ac^2)^{[n/2]}(ab^2)^{2\{n/2\}} .
\end{equation}
Here we set $d\star b=d\star c=a(cb+bc)$.
\end{prop}
\begin{proof}
We again proceed by induction on $n$. When $n=1$ the left hand
side of \eqref{equ:shuffle} is
 $$cd-c\sha d+d\star c
 =-ac(b+c)-a(b+c)c+abc+acb=-2ac^2.$$
Similarly
 $$bd-b\sha d+d\star b
 =-ab(b+c)-a(b+c)b+abc+acb=-2ab^2.$$
So the proposition is true when $n=1$. Now assume that
\eqref{equ:stuffle} holds up to $n-1$ for some $n\ge 2$. We will
use repeatedly the following recursive expression of the shuffle
product: for any letters $x,y$ and words $w_1$ and $w_2$:
\begin{equation}\label{equ:recursha}
 (xw_1)\sha(yw_2)=x\Big(w_1\sha(yw_2)\Big)+y\Big((xw_1)\sha
 w_2\Big).
\end{equation}
Thus
 \allowdisplaybreaks{
 \begin{align}
 \ &\sum_{i=0}^{2n} (-1)^i\sha_i\big((cd)^{\star n})\notag \\
 =&\sum_{j=0}^{n}
 (d\star c)^j\sha (cd)^{\star (n-j)}-\sum_{j=1}^{n}
 \Big(d\star (cd)^{(j-1)}\Big)\sha \Big(c(d\star c)^{n-j}\Big)\notag\\
 =&(d\star c)^n+(cd)^{\star n}+\sum_{j=1}^{n-1}
 \Big\{a\Big(\big((b+c)\star c(d\star c)^{j-1}\big)\sha (cd)^{\star (n-j)}\Big)
 + c\Big((d\star c)^j\sha \big((d\star
 c)^{n-j-1}d\big)\Big)\Big\}\notag\\
 -&\sum_{j=1}^{n}
  \Big\{ a\Big(\big((b+c)\star (cd)^{\star(j-1)}\big)\sha \big(c(d\star c)^{n-j}\big)\Big)
  + c\Big(\big(d\star (cd)^{(j-1)}\big)\sha  (d\star c)^{n-j}\Big) \Big\}\notag\\
=&a\sum_{j=1}^{n}
 \Big((bc+cb)(d\star c)^{j-1}\Big)\sha (cd)^{\star (n-j)} \notag\\
 \ &-a\Big( (b+c)\sha \big(c(d\star c)^{n-1}\big)\Big)
 -a\sum_{j=2}^{n}\Big((bc+cb)(d\star c)^{j-2}d\Big)\sha \Big(c(d\star c)^{n-j}\Big)\notag \\
 =&ab\sum_{j=1}^{n}
 \Big(c(d\star c)^{j-1}\Big)\sha (cd)^{\star (n-j)}
 +ac\sum_{j=1}^{n}\Big(b(d\star c)^{j-1}\Big)\sha (cd)^{\star (n-j)} \label{equ:ac1}\\
 +& ac\sum_{j=1}^{n-1} \Big((bc+cb)(d\star c)^{j-1}\Big)\sha \Big((d\star c)^{(n-j-1)}d\Big)  \label{equ:ac2}\\
 -&ab\sum_{j=1}^{n} \Big( (cd)^{\star(j-1)}\Big)\sha \Big(c(d\star c)^{n-j}\Big)
 -ac{\sum_{j=1}^{n}}\raisebox{0.3ex}{$'$} (bd)^{\star(j-1)}\sha \Big(c(d\star c)^{n-j}\Big) \label{equ:ac3}\\
 -&ac\Big( (b+c)\sha  (d\star c)^{n-1} \Big)
 -ac\sum_{j=2}^{n}\Big((bc+cb)(d\star c)^{j-2}d\Big)\sha
 \Big((d\star c)^{n-j}\Big) \label{equ:ac4}
 \end{align}}
 where in $\sum'$ we used the fact that $d\star c=d\star b$. Now cancelling the terms beginning with
 $ab$ and regrouping we get:
\begin{align}
\eqref{equ:ac1}+\eqref{equ:ac3}=& ac\sum_{j=1}^{n}
\Bigg\{\Big(b(d\star c)^{j-1}\Big)\sha (cd)^{\star
 (n-j)}-(bd)^{\star(j-1)}\sha \Big(c(d\star c)^{n-j}\Big) \Bigg\}\notag\\
=&acb\left\{\sum_{j=1}^{n}(d\star c)^{j-1} \sha (cd)^{\star
 (n-j)}-\sum_{j=2}^{n}\Big(d(cd)^{\star(j-2)}\Big)\sha \Big(c(d\star
 c)^{n-j}\Big)\right\} \label{acb1}\\
 +&ac^2\left\{\sum_{j=1}^{n-1} \Big(b(d\star c)^{j-1}\Big)\sha  \Big((d\star c)^{
 (n-j-1)}d\Big)-\sum_{j=1}^{n} (bd)^{\star(j-1)}\sha(d\star
 c)^{n-j}\right\} \label{ac^21}
\end{align}
Let us denote the right hand side of \eqref{equ:shuffle} by
$f_n(a,b,c)$. Notice that we can safely change $c$ to $b$ in the
second big bracket above and therefore by induction assumption we
get
 \begin{equation}\label{9+10}
\eqref{acb1}+\eqref{ac^21}=
 acb(f_{n-1}(a,b,c))-ac^2(f_{n-1}(a,c,b)) .
\end{equation}
Consider now the remaining terms in $\sum_{i=0}^{2n}
(-1)^i\sha_i\big((cd)^{\star n})$:
\begin{align*}
\eqref{equ:ac2}+\eqref{equ:ac4}=& ac\sum_{j=1}^{n-1} \Big((bc+cb)(d\star c)^{j-1}\Big)\sha \Big((d\star c)^{(n-j-1)}d\Big) \\
 -&ac\Big( (b+c)\sha  (d\star c)^{n-1} \Big)
 -ac\sum_{j=2}^{n}\Big((bc+cb)(d\star c)^{j-2}d\Big)\sha
 \Big((d\star c)^{n-j}\Big)
\end{align*}
By recursive formula \eqref{equ:recursha} the above expression can
be simplified to
\begin{align*}
    -&acb\left\{\sum_{j=1}^{n}(d\star c)^{j-1} \sha (cd)^{\star
 (n-j)}-\sum_{j=1}^{n-1} \Big(c(d\star c)^{j-1}\Big)\sha \Big((d\star c)^{(n-j-1)}d\Big)\right\} \\
 +&ac^2\left\{\sum_{j=1}^{n-1} \Big(b(d\star c)^{j-1}\Big)\sha  \Big((d\star c)^{
 (n-j-1)}d\Big)-\sum_{j=1}^{n}\Big((bd)^{\star(j-1)}\Big)\sha(d\star
 c)^{n-j}\right\} \\
 =&-acb(f_{n-1}(a,b,c))-ac^2(f_{n-1}(a,c,b)),
\end{align*}
where all the terms beginning with $aca$ are cancelled out. Adding
this to \eqref{9+10} we finally find
 $$\sum_{i=0}^{2n}
(-1)^i\sha_i\big((cd)^{\star n})=-2ac^2(f_{n-1}(a,c,b))=
 f_n(a,b,c).$$
This completes the proof of identity \eqref{equ:shuffle}. Notice
that throughout the above proof we may exchange $b$ and $c$ and
thus identity \eqref{equ:shuffleb} follows immediately. This
completes the proof of the proposition and therefore
Theorem~\ref{thm:iden}.
\end{proof}

If we consider the partial sums of Euler sums in
Theorem~\ref{thm:iden} then we get the following result due to
Bowein, Bradley and Broadhurst (see \cite[Conjecture 1]{BBB1}).
\begin{cor} \label{maincor}
Define a sequence $\{a_n(t)\}_{n\ge 1}$ by: $a_1(t)=a_2(t)=1$, and
recursively
\begin{equation}\label{equ:recur}
 n (n+1)^2 a_{n+2}= n (2n+1)a_{n+1}+(n^3+(-1)^{n+1}t)a_n, \quad\forall n\ge
 1.
\end{equation}
Then
\begin{equation}\label{equ:limitan}
 \lim_{n\to \infty}a_n(t)
 =\prod_{n=1}^\infty\left(1+\frac{t}{8n^3}\right).
\end{equation}
\end{cor}
\begin{proof}
It is easy to check that the sequence
$$
 \tilde{a}_n(t)=1+\sum_{i=1}^\infty t^i \sum_{n>l_1>k_1>\cdots>l_i>k_i\ge 1}
\frac{(-1)^{l_1+\cdots+l_i}}{l_1^2k_1\cdots l_i^2k_i}
$$
satisfies the initial conditions $\tilde{a}_1(t)=\tilde{a}_2(t)=1$
and the recursive relation  \eqref{equ:recur}. Hence
$a_n(t)=\tilde{a}_n(t)$ and
$$\lim_{n\to \infty}a_n(t)
 =1+\sum_{i=1}^\infty\zeta(\{\bar 2,1\}^i)t^i.$$
On the other hand, let
$$
b_n(t):=\prod_{i=1}^n\left(1+\frac{t}{8i^3}\right)=1+\sum_{i=1}^\infty
 \frac{t^i}{8^i} \sum_{n>l_1>\cdots>l_i\ge 1} \frac{1}{l_1^3\cdots l_i^3}.$$
Then \eqref{equ:limitan} is equivalent to
 \begin{equation}\label{equ:an=bn}
\lim_{n\to \infty}a_n(t)=\lim_{n\to \infty}b_n(t).
\end{equation}
But clearly
$$\lim_{n\to \infty}b_n(t)
 =1+\sum_{i=1}^\infty\zeta(\{3\}^i)\frac{t^i}{8^i}.$$
So \eqref{equ:an=bn} is equivalent to Theorem~\ref{thm:iden} and
the corollary follows.
\end{proof}

\noindent {\em Email:} zhaoj@eckerd.edu


\begin{thebibliography}{9}
\bibitem{B}
J.~M.~Borwein, Email to the author dated Oct. 08, 2006.

\bibitem{BBB1}
D.~Borwein, J.~M.~Borwein, and D.~M.~Bradley, Parametric Euler sum
identities, \textit{J. Math. Anal. Appl.}, \textbf{316} (2006),
328--338.

\bibitem{BoB}
D.~Bowman and D.~M.~Bradley, The algebra and combinatorics of
shuffles and multiple zeta values \emph{Journal of Combinatorial
Theory}, Series A, Vol. \textbf{97} (1)(2002), 43--61.

\bibitem{BdB}
J.~M.~Borwein, and D.~M.~Bradley, Thirty-Two Goldbach Variations,
\textit{International Journal of Number Theory},   \textbf{2}
(1)(2006), 65--103.

\bibitem{BBB2}
J.~M.~Borwein, D.~M.~Bradley and D.~J.~Broadhurst, ``Evaluations
of $k$-fold Euler/Zagier sums: A~compendium of results for
arbitrary $k$,'' {\em Electron.\ J.~Combin.} {\bf 4} (1997),
No.~2, \#R5.


\bibitem{EZface}
J.~Borwein, P.~Lisonek, and P.~Irvine, An interface for evaluation
of Euler sums, available online at
http://oldweb.cecm.sfu.ca/cgi-bin/EZFace/zetaform.cgi

\bibitem{Br1}
D. J. Broadhurst, Massive 3-loop Feynman diagrams reducible to SC*
primitives of algebras of the sixth root of unity, \emph{European
Phys. J. C (Fields)}  \textbf{8} (1999), 311--333

\bibitem{Br2}
\underline{\phantom{D. J. Broadhurst}},  Conjectured enumeration
of irreducible multiple zeta values, from knots and Feynman
diagrams, preprint hep-th9612012.

\bibitem{Chen}
K.-T.-Chen, Algebras of iterated path integrals and fundamental
groups, \textit{Trans. Amer. Math. Soc}.  \textbf{156} (1971),
359--379.

\bibitem{LE1} L.~Euler, Meditationes circa singulare serierum genus,
\textit{Novi.\ Comm.\ Acad.\ Sci.\ Petropolitanae}, \textbf{20}
(1775), 140--186.

\bibitem{LE2}
\underline{\phantom{L.~Euler}}, \textit{Opera Omnia}, ser.~1,
vol.~15, B.~G.~Teubner, Berlin, 1927.

\bibitem{Gcyc}
A.~B.~Goncharov, The double logarithm and Manin's complex for
modular curves, \emph{Math. Res. Letters} \textbf{4} (1997),
617--636.

\bibitem{MG}
A.~B.~Goncharov, Yu.~I.~Manin, Multiple zeta-motives and moduli
spaces $M_{0,n}$, \emph{Compositio Math.} \textbf{140} (2004),
1–-14.

\bibitem{H1} M.~E.~Hoffman, Multiple harmonic series,
\textit{Pacific J.~Math.}, \textbf{152} (2)(1992), 275--290.

\bibitem{H2}
\underline{\phantom{M.~E.~Hoffman}}, Quasi-shuffle products,
\textit{J. Algebraic Combin.} \textbf{11} (2000), 49--68.

\bibitem{H3}
\underline{\phantom{M.~E.~Hoffman}}, Algebra of Multiple Zeta
Values and Euler Sums, Mini-Conference on Zeta Functions, Index,
and Twisted K-Theory: Interactions with Physics, Oberwolfach,
Germany, May 2, 2006. Available through
www.usna.edu/Users/math/meh.

\bibitem{IKZ}
K. Ihara, M. Kaneko, and D. Zagier, Derivation and double shuffle
relations for multiple zeta values,  \textit{Compositio Math.}
\textbf{142} (2006), 307--338;

\bibitem{LM} T.~Q.~T.~Le and J.~Murakami, Kontsevich's integral for the Homfly
polynomial and relations between values of the multiple zeta
functions, \emph{Topology Appl.} \textbf{62} (1995), 193--206.

\bibitem{Zag} D.~Zagier, Values of Zeta
Function and Their Applications, \textit{Proceedings of the First
European Congress of Mathematics}, \textbf{2}, (1994), 497--512.

\bibitem{Zh} J.~Zhao, Double shuffle relations of special values of multiple
polylogarithms, arxiv:0707.1459

\bibitem{Zl} S.~Zlobin, Email to the author dated May 22, 2007.
\end{thebibliography}
\end{document}